\documentstyle[12pt]{article}

\newtheorem{theorem}{Theorem}[section]

\font\bbc=msbm10 scaled 1200

\begin{document}

\title{Asymptotics for the Almost Sure Lyapunov Exponent for the
Solution of the Parabolic Anderson Problem}

\author{Rene Carmona \\[1pt]
\normalsize Department of Operations Research and Financial Engineering\\[-4pt]
\normalsize Princeton University, Princeton, NJ 08544\\[7pt]
Leonid Koralov\footnote{Supported by the Institute for Advanced Study and
the NSF postdoctoral fellowship} \\[1pt]
\normalsize School of Mathematics \\[-4pt]
\normalsize IAS, Princeton, NJ 08540\\[7pt]
Stanislav Molchanov \\[1pt]
\normalsize Department of Mathematics \\[-4pt]
\normalsize UNC at Charlotte, Charlotte, NC 28223}

\date{}
\maketitle

\begin{abstract}
We consider the Stochastic Partial Differential Equation
\[
u_t =  \kappa \Delta u + \xi(t,x)u, ~~~ t
\geq 0, ~~~ x \in \mbox {\bbc Z}^d~.
\]
The potential is assumed to be Gaussian white noise in time, stationary in
space. We obtain the asymptotics of the almost sure
 Lyapunov exponent $\gamma(\kappa)$ for the solution as $\kappa \rightarrow
0$. Namely $\gamma(\kappa) \sim \frac{c_0}{\ln(1/\kappa)}$, where
the constant $c_0$ is determined by the correlation function of the
potential.
\medskip
\end{abstract}

\section {Introduction}
\label{se1}
    We study the long time asymptotic behavior for the solution of the
parabolic equation with random potential
\begin{equation}
 u_t =  \kappa \Delta u + \xi(t,x)u  \label{e1}
\end{equation}
\[
 t \geq 0, ~~~ x \in \mbox {\bbc Z}^d
\]
with  the initial data $u(0,x)$ concentrated at the point $x=0$.  The potential
$\xi(t,x)$ is  a mean zero stationary Gaussian field, which is white noise in
time. The almost sure Lyapunov exponent of the solution is defined as the
following limit
\begin{equation}
\gamma(\kappa) = \lim_{t \rightarrow \infty} \frac{\ln u(t,0)}{t}~, \label{x1}
\end{equation}
provided the limit exists a.s. and is nonrandom. The existence of the
Lyapunov exponent for the above type of potentials easily follows from
the Kingman's subadditive ergodic theorem \cite{St}.  Indeed, the fundamental
solution $q(s,x,t,y)$ satisfies
\[
q(s,0,u,0) \geq  q(s,0,t,0)q(t,0,u,0)
\]
when $s < t < u$, while
\[
\frac{{\rm E} \ln q(t,0,0,0)}{t} \leq \frac{ \ln  {\rm E} q(t,0,0,0)}{t}~,
\]
which is bounded uniformly in $t$ \cite{CM}. Thus, due to time stationarity of
the potential we can apply the subadditive ergodic theorem to the
process $\ln u(t,0)$, and thus prove the existence of the limit in (\ref{x1}). We are interested in the asymptotic behavior of the
Lyapunov exponent $\gamma(\kappa)$ for small $\kappa$. In \cite{CM} the
 bounds on   $\gamma(\kappa)$ were obtained with the
lower bound being of order $\frac{1}{\ln(1/\kappa)}$, and the upper bound of
order $\frac{\ln \ln(1/\kappa)}{\ln(1/\kappa)}$. In \cite{CVM} the upper bound
was improved to have the same order as the lower bound. In the present
 paper we
obtain the asymptotics for the Lyapunov exponent. That is, instead of the
upper and lower bounds with different constants,  we show that
a constant $c_0$ exists such that
\[
\gamma(\kappa) \sim \frac{c_0}{\ln(1/\kappa)}~~~{\rm as}~~~\kappa
 \rightarrow 0~.
\]
 We shall prove the following theorem:
\begin{theorem}
\label{tm1}
Let $ (W(t,x), t \geq 0, x \in \mbox {\bbc Z}^d )$ be a mean zero Gaussian
field with covariance ${\rm E} (W(t,x)W(s,y)) = \min(s,t) Q(x-y)$, where $Q(x)$ is
not identically constant. Then there
exists a constant $c_0 > 0 $ such that the
solution of the stochastic PDE
\begin{equation}
 d u(t,x) = \kappa \Delta u(t,x)dt + u(t,x) \circ dW(t,x)  \label{e2}
\end{equation}
with initial data $u(0,x) = \delta(x)$ satisfies almost surely
\[
\lim_{t \rightarrow \infty} \frac{\ln u(t,0)}{t}
= \gamma(\kappa) \sim \frac{c_0}{\ln 1/\kappa}~~~{\rm as}~~~\kappa \rightarrow
0~.
\]
\end{theorem}
The outline of the proof is the following:
Via Feynman-Kac formula we represent the solution $u(t,0)$ as an integral
 over the set of
paths of a continuous time random walk starting at $0$ and ending at $0$ at
 time $t$. The contribution from each path  is a functional of a
Gaussian random variable. We study a small subset of the set of
all paths which gives the main contribution to the solution.
In order to do that
we need to study the structure of the maxima of the Gaussian field defined over
the set of all paths of the random walk with a fixed number of jumps.

The paper is organized as follows.
In section \ref{se2} we introduce the notations. In section \ref{se3} we
prove the lower bound for the Lyapunov exponent, while in section
\ref{se4} we
prove the matching upper bound.
\section{Notations and Preliminary Considerations}
\label{se2}
 Let $P_m$ be the set of paths of a discrete time random walk of length $m$
starting and finishing at $0$. Let $p(t,m)$ be the probability
that at time $t$ the
Poisson process with intensity $2d\kappa$  is equal to $m$.
 Let $S(t,m)$ be the set of possible times of the
jumps: $\widetilde{t}= (t_1, ...,t_m)$ with $0 \leq t_1 \leq ... \leq t_m \leq
t$.
For $\widetilde{t} \in S(t,m)$ and
$\widetilde{x} = ( x_0, x_1, ..., x_m) \in P_m$ we define
 \begin{equation}
X_m( \widetilde{t},\widetilde{x}) =  W(t_1, x_0) - W (0,x_0) + W(t_2,x_1) -
W(t_1,x_1) +...+ W(t,x_m) - W(t_m, x_m)~. \label{y1}
\end{equation}
By the Feynman-Kac formula the solution of equation \ref{e2} can be
written as
\[
u(t,0) = \sum\limits_{m = 0}^{\infty} p(t,m) \sum\limits_{\widetilde{x} \in
P_m} \frac{1}{(2d)^m} \int\limits_{\widetilde{t} \in S(t,m)}
\frac{m !}{t^m} e^{X_m(\widetilde{t}, \widetilde{x})} dt_1 ... dt_m =
\]
\begin{equation}
\sum\limits_{m = 0}^{\infty} \kappa^m e^{-2d\kappa t}
 \sum\limits_{\widetilde{x} \in
P_m}  ~~\int\limits_{\widetilde{t} \in S(t,m)}  e^{X_m(\widetilde{t},
 \widetilde{x})} dt_1 ... dt_m~.
\label{e5}
\end{equation}
The proof of both lower and upper bounds will be a refinement of the
corresponding considerations in \cite{CM} and \cite{CVM}. The central idea
of the lower estimation in \cite{CM} is the following: with the fixed times
of the jumps $(t_1,..., t_m)$ one selects a path, which gives a significant
contribution to the integral in (\ref{e5}). The path is constructed as
follows. With $(x_0,...,x_{k-1})$ selected, one takes $x_k$ to be the
neighbor of $x_{k-1}$, which maximizes the increment $W(t_{k+1}, x_k) -
W(t_{k}, x_k)$. Such a path follows the local maximum of the potential.
Better estimations can be obtained by similar optimization with respect to
pairs, triples, etc. of random jumps. In this paper, rather than maximizing
the increments $W(t_{k+1}, x_k) -
W(t_{k}, x_k)$ at each time step, we shall consider the maximum of the filed
$ X_m(\widetilde{t},
 \widetilde{x})$ over all possible paths $\widetilde{x}$ and times of
the jumps $\widetilde{t}$. We further find the optimal number of jumps as a
function of $t$ and $\kappa$.

The proof of the upper estimate is similar to that of \cite{CVM}. However, now
we are able  to  keep track of the constant, so that to make sure that the
upper and lower bounds match, giving us the exact asymptotics for the
Lyapunov Exponent.

For $\alpha$ nonnegative let $m_0(\alpha,t) = [{\alpha t}]$ (the integer
part). Later $\alpha$ will be taken to be equal to $z {\ln^{-2}(1/\kappa)}$
 where
$z$ will vary. It is the paths with the  number of jumps of order
$  \ln^{-2}(1/\kappa) t$ that
give the main contribution to the positive solution $u(t,x)$. Let us stress
 that
for typical paths the number of jumps $m_1$ has the order $\kappa t$, that
is $m_0 >> m_1$.

 Let $T_{\alpha}$ be the space of paths of the continuous time random walk
starting and finishing at $x = 0$
with no
more than $m_0(\alpha, t)$ jumps. Thus to specify an element of $T_{\alpha}$ we need to
specify the number of jumps $m$,
 $ \widetilde{t} \in S(t,m)$, and $\widetilde{x} \in P_m$.
  Then $X_m(\widetilde{t},
\widetilde{x})$ is a Gaussian field over $T_{\alpha}$.
 Let
\begin{equation}
g(\alpha,t) =
\sup_{( m , \widetilde{t} , \widetilde{x}) \in  T_{\alpha}} X_m(\widetilde{t},
\widetilde{x})~, ~~{\rm and}~~ f(\alpha,t) = {\rm E} g(\alpha,t)~. \label{p1}
\end{equation}
Then $f(\alpha ,t)$ is a
 superadditive nonrandom function of $t$. In order to see
that $\lim_{t \rightarrow \infty} \frac{f(\alpha ,t)}{t} $ exists we  need to
show that $f$ does not grow faster than linearly in $t$.
 We shall use the following
entropy estimate \cite{Ad}:
For a mean zero
Gaussian field $X(\tau)$ over a set $T$ one defines the canonical
metric $d(\tau_1, \tau_2) = ({\rm E}(X(\tau_1) - X(\tau_2))^2)^{1/2}$.
 The metric
entropy $N(\varepsilon)$ is the smallest number of closed $d$-balls of radius
 $\varepsilon$
 needed to cover $T$. One has the following theorem \cite{Ad}
\begin{theorem}
\label{ent}
There exists a universal constant $K$ such that
\[
{\rm E} \sup_{T} X(\tau) \leq K \int\limits_{0}^{\infty}
\sqrt{\ln N(\varepsilon)} d \varepsilon
\]
provided that the RHS is finite.
\end{theorem}
Let us estimate the entropy function $N(\varepsilon)$ for the field $X_m(\widetilde{t},
\widetilde{x})$ defined over $T_{\alpha}$.
Below $c_1,$  $c_2,$ etc. denote constants, which depend only on the
 field $W(t,x)$, unless it is indicated otherwise. Note that
the diameter of $T_{\alpha}$ does not exceed $2\sqrt{t}$, since
${\rm E}(X_m(\widetilde{t},
\widetilde{x}))^2 = t$. Assume that $\varepsilon \leq 2\sqrt{t}$.
For $\widetilde{t} = (t_1,...,t_m)$ and $\widetilde{s} = (s_1,...,s_m)$, with
$m$ and $\widetilde{x} \in P_m$ fixed, the distance between
 $(m, \widetilde{t},
\widetilde{x})$ and $(m, \widetilde{s},
\widetilde{x})$ in the canonical metric is estimated as follows:
\begin{equation}
d((m, \widetilde{t}, \widetilde{x}),  (m, \widetilde{s},
\widetilde{x})) = \sqrt{ {\rm E}(X_m(\widetilde{t},
\widetilde{x}) - X_m(\widetilde{s},
\widetilde{x}) )^2 }  \leq c_1 \sqrt{\sum_{i = 1}^{m}|t_i - s_i|}~.
\label{est11}
\end{equation}
The  inequality is the statement of Lemma 2.1 of \cite{CVM}, it can also
be seen directly from (\ref{y1}). Consider the finite subset $U(t,m) $ of
$S(t,m)$ defined as follows: $\widetilde{t} = (t_1,...,t_m) \in U(t,m)$ if
$0 \leq t_1 \leq t_2 \leq ... \leq t_m \leq t$, and all $t_i$ are positive
integer multiples of
 $\frac{\varepsilon^2}{c_1^2 m} $.
By (\ref{est11}) the  elements $(m, \widetilde{t}, \widetilde{x})$
 with $\widetilde{t} \in U(t, m) $ form an $\varepsilon$-net
in $T_{\alpha}$. The number of elements in $U(t,m)$ is equal to
$\frac{([\frac{c_1^2 t m}{\varepsilon^2}]+m -1)!}{([\frac{c_1^2 t m}{\varepsilon^2}]-1)! m!}$, which, for $\varepsilon \leq 2\sqrt{t}$
 is estimated as follows
\[
\frac{([\frac{c_1^2 t m}{\varepsilon^2}]+m -1)!}{([\frac{c_1^2 t m}{\varepsilon^2}]-1)! m!} \leq
\frac{([\frac{c_1^2 t m}{\varepsilon^2}]+m -1)^m}{m!} \leq
(\frac{c_2 t m }{\varepsilon^2})^m \frac{1}{m!} \leq
(\frac{ c_3 t }{\varepsilon^2})^m~,
\]
where the last inequality is due to Stirling's formula. Since there are
less than
$(2d)^m$ elements in $P_m$ and $m \leq m_0$, the entropy is estimated as
follows for $\varepsilon \leq 2\sqrt{t}$
\[
N(\varepsilon) \leq
 \sum_{m = 0}^{m_0} (2d)^m (\frac{ c_3 t }{\varepsilon^2})^m
\leq   (\frac{ c_4 t }{\varepsilon^2})^{m_0} ~.
\]
Thus,
\begin{equation}
N( \varepsilon) \leq
\max
\{ 1, (\frac{ c_4 t }{\varepsilon^2})^{m_0} \}~~
{\rm for}~~{\rm all}~~ \varepsilon~.
\label{entrest}
\end{equation}
Therefore
\[
\int_0^{\infty} \sqrt{ \ln N(\varepsilon)} d \varepsilon \leq
\sqrt{m_0} \int_0^{\sqrt{c_4t}}
 \sqrt{\ln(\frac{c_4t}{\varepsilon^2})}d\varepsilon=
\sqrt{c_4 m_0 t} \int_{0}^{1} \sqrt{\ln(\frac{1}{\varepsilon^2})}d
\varepsilon  \leq c_5\sqrt{m_0t} ~.
\]
 The entropy estimate implies that
\[
f(\alpha  ,t) \leq c_6 \sqrt{m_0 t} =
c_6 \sqrt{[\alpha t] t}  \leq c_7(\alpha ) t~.
\]
Therefore the limit $\lim_{t \rightarrow \infty} \frac{f(\alpha ,t)}{t}$ exists and
is finite. It will be denoted by $F(\alpha )$.
Due to the fact that
$W(t,x)$ is a Wiener process in $t$, and is therefore self-similar,
\[
\lim_{t \rightarrow \infty} \frac{1}{t}{
{\rm E}  \sup_{T_{\alpha}} X_m (\widetilde{t}, \widetilde{x})
} =
\sqrt{\alpha}
\lim_{t \rightarrow \infty} \frac{1}{t}{{\rm E} \sup_{T_1} X_m (\widetilde{t}, \widetilde{x})~,
}
\]
which implies
\begin{equation}
F(\alpha ) = \sqrt{\alpha }  F(1)~. \label{frm1}
\end{equation}

In order to relate the maximum of a 'typical' realization of a Gaussian field
to the expectation of the maximum, controlled by the entropy estimate, one
uses the Borell's inequality \cite{Ad}
\begin{theorem}
\label{bor}
Let $X(\tau)$ be a centered Gaussian process with sample paths bounded almost
surely. Then ${\rm E} \sup_{T} X(\tau) < \infty$,  and for all $\lambda  > 0$
\begin{equation}
{\rm Prob} \{ | \sup_{T} X(\tau) - {\rm E} \sup_{T} X(\tau)| > \lambda \} \leq
2 e^{-\frac{1}{2} \lambda^2/ \sigma^2}~, \label{be}
\end{equation}
where $\sigma^2 = \sup_{T} {\rm E} X^2(\tau)$.
\end{theorem}
We next make an intuitive argument providing us with the
asymptotics for the Lyapunov exponent. The rigorous derivation will be given
in the following two sections. We claim that the main contribution to the
sum in the RHS of (\ref{e5}) comes from the terms
 with $m \sim z{\ln^{-2}(1/\kappa)}t$,
where $z$ does not depend on $\kappa$.
Moreover, since we are taking the logarithm of $u$, only the
factors $\kappa^m$ and $e^{X_m(\widetilde{t}, \widetilde{x})}$ are important.
We substitute the latter by $e^{F(z{\ln^{-2}(1/\kappa)}) t }$.
Taking the logarithm of the product of these two factors, that is ignoring
the summation and integration in the RHS of (\ref{e5}), we arrive at the
expression
\[
\ln (\kappa^m  e^{F(z{\ln^{-2}(1/\kappa)})t })~.
\]
Dividing by $t$ we obtain
\[
\ln(\kappa^{\frac{z}{\ln^2 (1/\kappa)}}
e^{\frac{\sqrt{z}}{\ln (1/\kappa)}F(1)}
) = \frac{F(1)\sqrt{z} -z}{\ln(1/\kappa)}~.
\]
The maximum of this function is equal to $\frac{F(1)^2}{4\ln(1/\kappa)}$ and
is achieved at $z = F(1)^2/4$. We shall prove that $\gamma(\kappa) \sim
\frac{F(1)^2}{4\ln(1/\kappa)}$. Note that $F(1)$ is greater than zero if
$Q(x)$ is not identically constant.

\section{Proof of the Lower Bound}
\label{se3}
We select $c_0$ (which is the same constant as in the statement of theorem
\ref{tm1}) as follows: $c_0 = \frac{F(1)^2}{4}$.
Let $\delta_0 > 0$ be given.
In order to prove the estimate from below it is
sufficient to show that there exists $\kappa_0$ such that
\begin{equation}
\gamma(\kappa) \geq \frac{c_0}{\ln(1/\kappa)}(1- \delta_0)~~{\rm for}~~
\kappa \leq \kappa_0~. \label{fr1}
\end{equation}
Let $\alpha (\kappa) = \frac{c_0}{\ln^2 (1/\kappa)}$.
Since we know that the Lyapunov exponent exists and is nonrandom in order to
prove (\ref{fr1}) it is
sufficient to show that for $\kappa \leq \kappa_0$, for large enough $t$
\begin{equation}
{\rm Prob} \{
\frac{1}{t} \ln ( \sum_{m \leq [ \alpha (\kappa)  t]}
\kappa^m e^{-2d\kappa t} \sum_{\widetilde{x} \in P_m}~~
\int\limits_{\widetilde{t} \in S(t,m)} e^{X_m(\widetilde{t},\widetilde{x})}
dt_1...dt_m ) \geq \frac{c_0(1-\delta_0)}{\ln(1/\kappa)} \}
\geq 1/2~.~ \label{t1}
\end{equation}
Since $\frac{m \ln \kappa}{t} \geq \frac{-c_0}{\ln(1/\kappa)}$
for $m \leq [\alpha(\kappa) t]$,
the LHS of the first inequality in (\ref{t1}) can be estimated from below
by the following expression
\[
- \frac{c_0}{\ln(1/\kappa)}  - 2d\kappa + \frac{1}{t}
\ln(\sup_{m \leq [ \alpha (\kappa) t ]}
\sup_{\widetilde{x} \in P_m} ~~\int\limits_{\widetilde{t} \in S(t,m)}
e^{X_m(\widetilde{t},\widetilde{x})}
dt_1...dt_m )~.
\]
Select $\kappa_1$ small enough so that
\begin{equation}
2d\kappa < \frac{c_0 \delta_0}
{2 \ln(1/ \kappa)} ~~{\rm for} ~~\kappa \leq \kappa_1~. \label{frm3}
\end{equation}
Thus what we want to show is that for sufficiently small $\kappa$
\begin{equation}
{\rm Prob} \{
 \frac{1}{t}
\ln(\sup_{m \leq [ \alpha(\kappa) t ]}
\sup_{\widetilde{x} \in P_m} ~~\int\limits_{\widetilde{t} \in S(t,m)} e^{X_m(\widetilde{t},\widetilde{x})}
dt_1...dt_m ) \geq \frac{c_0}{\ln(1/\kappa)}(2-\frac{\delta_0}{2}) \}
\geq 1/2
\label{est1}
\end{equation}
for large $t$. \par
We first consider the field $X_m(\widetilde{t},\widetilde{x})$ to be defined
on $T_{\alpha}^r$, a space of paths somewhat smaller than $T_{\alpha}$.
Namely $T^r_{\alpha}$ is the space of paths of the continuous time random
walk starting and finishing at $x=0$ with no more than $m_0(\alpha,t)
 = [\alpha t]$
jumps, with the jumps separated from each other and from
the endpoints of the interval $[0,t]$ by a distance of at least $2r$. Let
\[
g^r(\alpha, t) = \sup_{T_{\alpha}^r} X_m(\widetilde{t},\widetilde{x})~~
{\rm and}~~ f^r(\alpha, t) = {\rm E} g^r(\alpha, t)~.
\]
Then $~f^r(\alpha, t)~$ is a superadditive nonrandom function of $~t$ and there
exists the limit
$F^r(\alpha) =
\lim_{t \rightarrow \infty} \frac{f^r(\alpha, t)}{t}$. As in (\ref{frm1})
\[
F^r(\alpha) = \sqrt{\alpha} F^{\alpha r} (1) ~.
\]
We next demonstrate that $F^r(1) \rightarrow F(1)$ as $r \rightarrow 0$.
Indeed, let an arbitrary $\varepsilon >0 $ be given. Select $t_0$ such that
$F(1) - \frac{f(1, t_0)}{t_0} < \varepsilon/2$. Since
$f^r(1, t_0) \rightarrow f(1, t_0)$ for fixed $t_0$, we can find
$r_0$ such that $F(1) - \frac{f^r(1, t_0)}{t_0} < \varepsilon$ for
$r \leq r_0$. Note that $\frac{f^r(1, t)}{t}$ is increasing in $t$.
Therefore $F(1) - \lim_{t \rightarrow \infty} \frac{f^r(1, t)}{t} <
\varepsilon$ for $r \leq r_0$. Thus $F(1) -F^r(1) < \varepsilon$ for $r \leq
r_0$, which proves our claim.

Now we estimate from below the supremum of the field $ X_m(\widetilde{t},
\widetilde{x}) $ over $T^1_{\alpha}$. \\
Let
$\delta_1 = \frac{\delta_0 \sqrt{c_0}}{8}$. Select $\kappa_2$ such that
\begin{equation}
F^{\alpha(\kappa)}(1) \geq F(1) - \frac{\delta_1}{2} ~~{\rm for}~~ \kappa
\leq \kappa_2~. \label{frm3a}
\end{equation}
Then for $\kappa \leq \kappa_2$
there exists $ t_1 = t_1(\delta_1, \kappa)$, such that for $t \geq t_1$
\[
{\rm E} \sup_{T^1_{\alpha(\kappa)}}  X_m(\widetilde{t},
\widetilde{x}) \geq t(F^1(\alpha(\kappa)) -
\frac{\sqrt{\alpha(\kappa)} \delta_1}{2}) =
\]
\[
t(\sqrt{\alpha(\kappa)} F^{\alpha(\kappa)
}(1) - \frac{\sqrt{\alpha(\kappa)} \delta_1}{2}) \geq
\sqrt{\alpha(\kappa)} t (F(1) - \delta_1)~.
\]
Then by Borell's inequality (\ref{be}) with $\sigma^2 = t$ and $\lambda =
\sqrt{\alpha(\kappa)} \delta_1 t $
\begin{equation}
{\rm Prob} \{ \sup_{T^1_{\alpha(\kappa)} }  X_m(\widetilde{t},
\widetilde{x}) < \sqrt{\alpha(\kappa)} t (F(1) - 2 \delta_1) \} \leq 2e^{
\frac{-\delta_1^2 \alpha t}{2} }\leq \frac{1}{4}~, \label{eqe1}
\end{equation}
where the last inequality holds for $t$ sufficiently large.
The fact that we can estimate the supremum  of the field $ X_m(\widetilde{t},
\widetilde{x})$  from below  does not immediately allow us to prove
(\ref{est1}). The problem is that the expression in the LHS of the first
inequality in
 (\ref{est1}) involves integration in $t_1,...,t_m$, and thus we need to study
the fluctuations  of the field $ X_m(\widetilde{t},
\widetilde{x})$ as $\widetilde{t}$ is varied.
Consider the field
\[
Y(\widetilde{t_1}, \widetilde{t_2}, m , \widetilde{x}) =
X_m(\widetilde{t_1},
\widetilde{x}) - X_m(\widetilde{t_2},
\widetilde{x})~
\]
defined on the set $U_{\alpha, \beta} \subset T_{\alpha} \times T_{\alpha}$.
We say that $((\widetilde{t_1}, m_1, \widetilde{x_1}),
(\widetilde{t_2}, m_2, \widetilde{x_2})) \in U_{\alpha , \beta}$ if and only if
$m_1 = m_2 \leq [\alpha t],  \widetilde{x_1} = \widetilde{x_2}$, and
$|\widetilde{t^i_1} - \widetilde{t^i_2}| \leq \beta $ for $1 \leq i \leq m_1$.
Let
\[
G(\alpha, \beta) = \limsup_{t \rightarrow \infty} \frac{1}{t}  {\rm E} \sup_{U_{\alpha,
\beta}} Y~.
\]
 Then $G(\alpha, \beta) = \sqrt{\alpha} G(1,\alpha \beta) $. Below we
are going to show that $G(1, \beta) \rightarrow 0$ as $\beta \rightarrow 0$.
First, assuming that
 we have this statement, we prove (\ref{fr1}).
We take $\kappa_3$ so small that
\begin{equation}
G(1, \alpha(\kappa)) < \delta_1 ~~{\rm  for}~~ \kappa \leq \kappa_3~.
\label{frm3b}
\end{equation}
 Then for
$t$ sufficiently large ${\rm E} \sup_{U_{\alpha(\kappa) , 1}} Y \leq \delta_1
\sqrt{\alpha(\kappa)} t$. Thus by Borell's inequality
 with $\sigma^2 \leq  2t$ and $\lambda =
\sqrt{\alpha(\kappa)} \delta_1 t$
\begin{equation}
{\rm Prob} \{\sup_{U_{\alpha(\kappa)
 , 1}}Y  > 2 \sqrt{\alpha(\kappa)} t \delta_1 \} \leq
2 e^{\frac{-\delta_1^2 \alpha(\kappa) t}{4}} \leq \frac{1}{4} \label{aa1}
\end{equation}
where the last inequality holds for $t$ sufficiently large.
Take $\kappa_0 = \min \{ \kappa_1, \kappa_2, \kappa_3 \}$, so that
(\ref{frm3}), (\ref{frm3a}), and (\ref{frm3b}) hold.

From (\ref{aa1}) and (\ref{eqe1}) it follows that for
$\kappa \leq \kappa_0$  for $t$ sufficiently large  with probability of
at least
 $1/2$ one has
\[
\sup_{T^1_{\alpha(\kappa)}} X > \sqrt{\alpha(\kappa)} t(F(1) -
 2 \delta_1) ~~{\rm and}~{\rm at}~ {\rm the}~
{\rm same}~{\rm   time}~~\sup_{U_{\alpha(\kappa), 1}}Y
 < 2 \sqrt{\alpha(\kappa)} t \delta_1~.
\]
Therefore for large $t$ on a set of probability of at least
$1/2$
\begin{equation}
\frac{1}{t}
\ln(\sup_{m \leq [ \alpha(\kappa) t ]}
\sup_{\widetilde{x} \in P_m} ~~\int\limits_{\widetilde{t} \in S(t,m)}
 e^{X_m(\widetilde{t},\widetilde{x})}
dt_1...dt_m ) \geq
\frac{1}{t}  \ln (
e^{t \sqrt{\alpha(\kappa)}(F(1) - 4 \delta_1)}) =
\sqrt{\alpha(\kappa)}(F(1) -4 \delta_1). \label{eq6}
\end{equation}
Recall that $F(1) = 2 \sqrt{c_0}$,  $\sqrt{\alpha(\kappa)} = \frac{\sqrt{c_0}}{
\ln(1/ \kappa)}$, and $\delta_1 = \delta_0 \sqrt{c_0}
 /8$. Thus
the RHS of (\ref{eq6}) is greater or equal than $\frac{c_0}{\ln (1/ \kappa)}
(2 - \frac{\delta_0}{2})$. Thus the
estimate (\ref{est1}) holds. It remains to show that $G(1, \beta) \rightarrow
0$ as $\beta \rightarrow 0$.

By (\ref{entrest}) the entropy function $N_X(\varepsilon)$ for the field
$X_m(\widetilde{t}, \widetilde{x})$ defined over $T_1$ (that is
with $m_0 = [t]$) is  estimated as
$N_X(\varepsilon) \leq \max \{ 1, (\frac{ct}{\varepsilon^2})^t \} $. The entropy function
$N_Y(\varepsilon)$ for the field $Y(\widetilde{t_1},  \widetilde{t_2}, m,
\widetilde{x})$ defined on the set $U_{1, \beta}$ is estimated as
$N_Y(\varepsilon) \leq N_X^2(\varepsilon / 2)$. The diameter of $U_{1, \beta}$
in the canonical metric associated with $Y$ does not exceed
$2 \sup_{U_{1, \beta}} ({\rm E} Y^2)^{1/2} \leq 2 \sqrt{ \beta t}$. Thus for
$\beta < c$ by the
entropy estimate
\[
{\rm E} \sup_{U_{1, \beta}} Y \leq
 K \int_0^{2 \sqrt{\beta t}} \sqrt{ \ln N_Y(\varepsilon)} d \varepsilon \leq
K \int_0^{2 \sqrt{\beta t}} \sqrt{ \ln (\frac{4 ct}{\varepsilon^2})^{2t} }
 d \varepsilon =
2 \sqrt{2 c} K t \int_0^{\sqrt{\beta/c }} \sqrt{\ln(\frac{1}{\varepsilon^2})}
d \varepsilon~.
\]
Then,
\[
G(1, \beta) = \limsup_{t \rightarrow \infty} \frac{1}{t} {\rm E} \sup_{U_{1, \beta}}
Y \leq 2 \sqrt{2 c} K \int_0^{\sqrt{\beta/c}}
 \sqrt{\ln(\frac{1}{\varepsilon^2})}
d \varepsilon~,
\]
which tends to $0$ as $\beta \rightarrow 0$.
 This completes the proof of (\ref{fr1}).
\section{Proof of the Upper Bound}
\label{se4}
Let $\delta_0 > 0$ be given.
In order to prove the estimate from above it is
sufficient to show that there exists $\kappa_0$ such that
\begin{equation}
\gamma(\kappa) \leq \frac{c_0}{\ln(1/\kappa)}(1+ \delta_0)~~{\rm for}~~
\kappa \leq \kappa_0~. \label{fr11}
\end{equation}
As in section \ref{se3} we ust the Feynman-Kac formula for the solution.
Since we know that the Lyapunov exponent exists and is nonrandom in order to
prove (\ref{fr11}) it is
sufficient to show that for $\kappa \leq \kappa_0$ for large enough $t$
\begin{equation}
{\rm Prob} \{
\frac{1}{t} \ln ( \sum_{m =0}^{\infty}
\kappa^m e^{-2d\kappa t} \sum_{\widetilde{x} \in P_m}~~
\int\limits_{\widetilde{t} \in S(t,m)} e^{X_m(\widetilde{t},\widetilde{x})}
dt_1...dt_m ) \leq \frac{c_0(1+\delta_0)}{\ln(1/\kappa)} \} \geq 1/2~.\label{xx}
\end{equation}
The $m$ summation in the LHS of the first
inequality of (\ref{xx})  will be performed over disjoint
intervals separately. For $\gamma >0$ let us estimate the probability of the
following event
\begin{equation}
A_{a_1,a_2}(\gamma)=
\{\sum_{[\frac{a_1 t}{\ln^2(1/\kappa)}
 +1] \leq m \leq [\frac{a_2 t}{\ln^2(1/\kappa)}]} \kappa^m e^{-2d\kappa t}
\sum_{\widetilde{x} \in P_m}~~
\int\limits_{\widetilde{t} \in S(t,m)} e^{X_m(\widetilde{t},\widetilde{x})}
dt_1...dt_m > e^{\gamma t} \} \label{in1}
\end{equation}
Calculating the number of terms in each of the sums and the
 area of the domain of
integration, we estimate the LHS  of the inequality in
 (\ref{in1}) from above by
\begin{equation}
t(a_2-a_1) \ln^{-2}(1/\kappa) \kappa^{a_1 t} e^{-2d \kappa t} (2d)^{
\frac{a_2 t}{\ln^2(1/\kappa)} }
(\sup_{m \leq [\frac{a_2t}{\ln^2(1/\kappa)}]} \frac{t^m}{m!})
\exp({\sup_{m \leq [\frac{a_2 t}{\ln^2(1/\kappa)}] } X_m})
\label{er}
\end{equation}
Since $t(a_2 - a_1) \ln^{-2}(1/\kappa) \leq (2d)^{
\frac{a_2 t}{\ln^2(1/\kappa)} }$
the logarithm of (\ref{er}) is estimated from above by
\begin{equation}
\frac{2a_2 t}{\ln^2(1/\kappa)} \ln(2d) -
a_1t \ln(1/\kappa) -2d \kappa t  +
\ln \sup_{m \leq [\frac{a_2 t}{\ln^2(1/\kappa)}]}
(\frac{t^m}{m!}) + \sup_{m \leq [\frac{a_2 t}{\ln^2(1/\kappa)}]}
 X_m(\widetilde{t}, \widetilde{x}), \label{fr15}
\end{equation}
Let $c_1$, $c_2$, etc. denote constants which may depend only on the
dimension $d$.  From Stirling's formula it follows that
\[
\ln \sup_{m \leq [\frac{a_2 t}{\ln^2(1/\kappa)}]}
(\frac{t^m}{m!}) \leq \sup_{m \leq [\frac{a_2 t}{\ln^2(1/\kappa)}]}
m(\ln t - \ln m + c_1) \leq \frac{c_2 \max(1, a_2) t
 \ln(\ln^2(1/\kappa))}{\ln^2(1/\kappa)}~.
\]
Therefore the quantity in (\ref{fr15}) is estimated from above by
\[
\sup_{m \leq [\frac{a_2 t}{\ln^2(1/\kappa)}]
} X_m(\widetilde{t}, \widetilde{x})-
a_1t \ln(1/\kappa)+ \frac{c_3 \max(1, a_2) t
 \ln(\ln^2(1/\kappa))}{\ln^2(1/\kappa)}~,
\]
 Thus the probability of the event
$A_{a_1, a_2}(\gamma)$ is estimated from above by
\[
{\rm Prob}  \{ \sup_{m \leq [\frac{a_2 t}{\ln^2(1/\kappa)}]} X_m(\widetilde{t},
\widetilde{x})-
a_1t \ln(1/\kappa)+ \frac{c_3 \max(1, a_2) t
 \ln(\ln^2(1/\kappa))}{\ln^2(1/\kappa)} > \gamma t \} =
\]
\begin{equation}
{\rm Prob}  \{ \sup_{m \leq [\frac{a_2 t}{\ln^2(1/\kappa)}]} X_m(\widetilde{t},
\widetilde{x})>
t(\gamma + a_1 \ln(1/\kappa) -
 \frac{o(\kappa) \max (1, a_2)}{\ln(1/\kappa)} ) \}~.
\label{in2}
\end{equation}
Recall that ${\rm E} \sup_{m \leq [\frac{a_2 t}{\ln^2(1/\kappa)}]}
 X_m(\widetilde{t}, \widetilde{x})$
is a superadditive function of $t$ with the limit
\[
\lim_{t \rightarrow \infty}
 \frac{
{\rm E}
  \sup_{m \leq [\frac{a_2 t}{\ln^2(1/\kappa)}]}
  X_m(\widetilde{t},
\widetilde{x})
}
{t}
 =
\frac{\sqrt{a_2}}{\ln(1/\kappa)} F(1)~.
\]
Thus ${\rm E}\sup_{m \leq [\frac{a_2 t}{\ln^2(1/\kappa)}]} X_m(\widetilde{t} ,
\widetilde{x})
\leq  \frac{t\sqrt{a_2}}{\ln(1/\kappa)} F(1)$. Therefore by Borell's inequality
(\ref{be}) with $\sigma^2 = t$
the probability in (\ref{in2}) is estimated from above by
\begin{equation}
2 \exp(-\frac{1}{2} t(\gamma + a_1 \ln(1/\kappa) -
 \frac{\sqrt{a_2}}{\ln(1/\kappa)} F(1) -
\frac{o(\kappa)\max(1, a_2)}{\ln(1/\kappa)})^2 )
\label{ex1}
\end{equation}
Take
\[
\gamma = \frac{c_0(1+\delta_0)}{\ln(1/\kappa)},~~~{\rm and}~~~
\gamma_n = \frac{c_0(1+\frac{\delta_0}{2}) - n \varepsilon}{\ln(1/\kappa)},
\]
where $\varepsilon  $ is a positive number to be selected below. We
 cover the axis $[0, \infty)$ by
intervals
\[
[a_1^n, a_2^n] = [(n-1)\varepsilon_1, n\varepsilon_1], ~~n \geq 1~,
\]
where $\varepsilon_1 $ is to be specified below.
Then for large $t$
\begin{equation}
P( \sum_{m =0}^{\infty}
\kappa^m e^{-2d\kappa t} \sum_{\widetilde{x} \in P_m}~~
\int\limits_{\widetilde{t} \in S(t,m)} e^{X_m(\widetilde{t},\widetilde{x})}
dt_1...dt_m > e^{\gamma t}) \leq
\sum_{n=1}^{\infty} P(A_{a_1^n, a_2^n}(\gamma_n)),
\label{in3}
\end{equation}
since
\[
\sum_{n=1}^{\infty} e^{\gamma_n t} < e^{\gamma t} ~~~{\rm for}~~{\rm large}~~
t~.
\]
Each term in the RHS of (\ref{in3}) is estimated by the expression of the
form (\ref{ex1}). Thus, in order to demonstrate that (\ref{xx}) holds it
is enough to show that
\begin{equation}
\sum_{n=1}^{\infty} 2 \exp(-\frac{1}{2} t(\gamma^n +
\frac{ a_1^n}{\ln(1/\kappa)} -
 \frac{\sqrt{a_2^n}}{\ln(1/\kappa)} F(1) -
\frac{o(\kappa)\max(1, a_2^n)}{\ln(1/\kappa)})^2 ) < \frac{1}{2}
\label{lt}
\end{equation}
for large $t$.  Recalling the definition of $ a^n_1, a^n_2,$ and $ \gamma^n$
we see that
\[
\gamma^n + \frac{ a_1^n}{ \ln(1/\kappa)} -
 \frac{\sqrt{a_2^n}}{\ln(1/\kappa)} F(1) -
\frac{o(\kappa)\max(1, a_2^n)}{\ln(1/\kappa)}=
\]
\begin{equation}
\frac{1}{\ln(1/\kappa)} (
c_0(1 + \frac{\delta_0}{2}) - n \varepsilon + n\varepsilon_1 - 2
\sqrt{n \varepsilon_1 c_0} - \varepsilon_1 - o(\kappa)
 \max(1 ,n \varepsilon_1))~.
\label{lt1}
\end{equation}
Since $c_0 +z - 2\sqrt{c_0 z} \geq 0 $ for all $z$ the expression in
(\ref{lt1}) is estimated from below by
\[
\frac{1}{\ln(1/\kappa)} (
c_0\frac{\delta_0}{2} - n \varepsilon + n \varepsilon_1 - \varepsilon_1
-o(\kappa) \max(1, n\varepsilon_1))~.
\]
Take such $\kappa_0$ that $o(\kappa) < \min(\frac{c_0\delta_0}{4},
 \frac{1}{3}) $ for $\kappa \leq \kappa_0$. Take $\varepsilon_1 =
 \frac{c_0 \delta_0}{4}$ and $\varepsilon = \frac{\varepsilon_1}{3}$. Then the
last expression is estimated from below by $
\frac{1}{\ln(1/\kappa)} n \frac{c_0 \delta_0}{12}$.
Therefore the LHS of (\ref{lt}) is estimated from above by
\[
\sum_{n=1}^{\infty} 2 \exp(-\frac{t}{2 \ln^2(1/\kappa)}
 (n\frac{c_0 \delta_0}{12})^2)~,
\]
which can be made arbitrarily small by selecting $t$ large enough.
This completes the proof of the upper bound.

\end{document}